\documentclass[12pt]{article}
\usepackage[utf8]{inputenc} 
\usepackage{amsmath, amsfonts, amssymb, amsxtra, amsthm}
\usepackage{lipsum}
\usepackage{graphicx} 
\usepackage{booktabs} 
\usepackage{array} 
\usepackage{paralist} 
\usepackage{verbatim} 
\usepackage{subfig}
\usepackage[a4paper]{geometry}
\newtheoremstyle{mystyle}{}{}{\itshape}{}{\bfseries}{.}{ }{\thmname{#1}\thmnumber{ #2}\thmnote{ (\bfseries #3)}}
\theoremstyle{mystyle}
\newtheorem{thm}{Theorem}

\newtheorem{lm}{Lemma}
\newtheorem{prop}{Proposition}
\theoremstyle{definition}
\newtheorem*{defi}{Definition}
\newtheorem*{prob}{Problem}
\theoremstyle{remark}
\newtheorem*{rmk}{Remark}

\DeclareMathOperator{\re}{Re}

\begin{document}        
\title{{\bf Smale's mean value conjecture for finite Blaschke products}}

\author{
Tuen-Wai Ng\thanks{
Partially supported by RGC grant HKU 704611P and HKU 703313P.
\mbox{\hspace{4truecm}}}
 \, and Yongquan Zhang\thanks{
Partially supported by a summer research fellowship of Faculty of Science, HKU.
\mbox{\hspace{0.9truecm}}
}}

\date{August 26, 2016}
\maketitle

\begin{center}
{\large\it Dedicated to David Minda on the occasion of his retirement}
\end{center}

\begin{figure}[b]

\rule[-2.5truemm]{5cm}{0.1truemm}\\[2mm]
{\footnotesize
2010{\it Mathematics Subject Classification: Primary} 30C15, 30D50.
\par\noindent
Department of Mathematics,
The University of Hong Kong,
Pokfulam Road.
\smallskip
\par\noindent E-mail: ntw@maths.hku.hk, zyqbsc@connect.hku.hk.
\par
}

\end{figure}

\begin{abstract}
Motivated by a dictionary between polynomials and finite Blaschke products, we study both Smale's mean value conjecture and its dual conjecture for finite Blaschke products in this paper. Our result on  the dual conjecture for finite Blaschke products allows us to improve a bound obtained by V. Dubinin and T. Sugawa for the dual mean value conjecture for polynomials. 
\end{abstract}



\section{Introduction}

\indent
In Stephen Smale's seminal paper \cite{smale} on the efficiency of Newton's method for approximating zeros of a non-constant complex polynomial $P$ of one complex variable, Smale proposed to study the smallest positive constant $c$ such that for any point $a$ in the complex plane $\mathbb{C}$, $\left|\frac{P(b)-P(a)}{b-a}\right| \le c|P'(a)|$ for at least one critical point $b$ of $P$ (zero of $P'$).  Let $M$ be  the least possible value of the factor $c$ for all non-linear polynomials and $M_n$ be the corresponding value for polynomials of degree $n$. It was proven by Smale \cite{smale} that $1 \le M\le4$ and he conjectured that $M=1$ or even $M_n=\frac{n-1}{n}$ and pointed out that the number $n-1 \over n$ would, if true, be the best possible bound here as it is attained (for any nonzero $\lambda$) when $P(z) = z^n-\lambda z$ and $a=0$. The conjecture was repeated in \cite{S85,SS86} and it is also listed as one of the three minor problems in Smale's famous problem list \cite{S00}. The conjecture is now known as Smale's mean value conjecture which has remained open since $1981$ even though it was proven to be true for many classes of polynomials (see \cite{CorRu},\cite{GS02},\cite{HiKa1},\cite{HiKa2},\cite{MP12},\cite{RS02},\cite{sheilsmall},\cite{Tis89},\cite{JT05} and \cite{WangY}).

There have been a number of refinements on the bound for $M$ over the past years. Beardon, Minda, and Ng \cite{BeardonMindaNg} showed that $M_n \le4^{(n-2)/(n-1)}$ using the hyperbolic metric on certain domains. Conte, Fujikawa and Lakic \cite{conteFujikawaLakic} proved that $M_n \le4\frac{n-1}{n+1}$ by referring to Bieberbach's coefficient inequality $|a_2| \le 2$ for univalent functions in the unit disk. Combining these two methods, Fujikawa and Sugawa \cite{fujikawaSugawa} verified that $M_n \le4\frac{1+(n-2)4^{1/(n-1)}}{n+1}$. Asymptotically, all these upper bounds of $M$ are of the form $4 - O({1 \over n})$ as $n \to \infty$. At present the best known result for large $n$ was obtained by Crane \cite{C07} who showed that  for $n \ge 8$, $M_n < 4 - {2.263 \over \sqrt{n}}$. All these upper bounds are close to $4$ when $n$ is large. On the hand, it was proved in \cite{Ng03} that one can replace $4$ by $2$ for a large class of polynomials which includes all odd polynomials with zero constant term. 

Notice that it is easy to show that Smale's mean value conjecture is equivalent to the following normalized conjecture:

\noindent
Let $P$ be a  monic polynomial of degree $n \ge 2$ such that $P(0)=0$ and $P'(0)=1$. Let $b_1, \dots , b_{n-1}$ be its critical points. Then

$$\min_{i}\left|{P(b_i) \over b_iP'(0)}\right| \le {n-1 \over n} \, . $$\

The following dual mean value conjecture was considered independently by the first author \cite{Ng07} and Dubinin and Sugawa \cite{dubininSugawa} around the same time. 

\noindent
{\bf Dual Mean Value Conjecture:} Let $P$ be a  monic polynomial of degree $n \ge 2$ such that $P(0)=0$ and $P'(0)=1$. Let $b_1, \dots , b_{n-1}$ be its critical points. Then

$$\max_{i}\left|{P(b_i) \over b_iP'(0)}\right| \ge {1 \over n} \, . $$\

Let the largest lower bound for $\max_{i}\left|{P(b_i) \over b_iP'(0)}\right|$ be $N$  when we  consider all polynomials satisfying $P(0)=0$ and $P'(0)=1$ and $N_n$ when we further restrict the degree of the polynomials to be $n$. By applying the theory of amoeba, the first author was able to show that $N_n>0$ (\cite{Ng07}) while Dubinin and Sugawa \cite{dubininSugawa} were able to show that $N_n \ge \frac{1}{n4^n}$. In this paper, we will study similar conjectures for finite Blaschke products.

It was first noted by Walsh \cite{walsh} that finite Blaschke products can be viewed as {\it non-euclidean polynomials} in the standard unit disk $\mathbb{D}$ and he has proven a version of Gauss-Lucas Theorem for finite Blaschke products. This point of view was also propagated by Beardon and Minda in \cite{BeardonMinda}, as well as Singer in \cite{Singer06}. Recently, a dictionary between polynomials and finite Blaschke products has been established by Ng and Tsang in \cite{NgTsang} based on the papers \cite{NgWang}, \cite{NgTsang15}, \cite{NgTsangZ} and see also \cite{CheungNgYam} and \cite{WangM} for more recent results not included in this dictionary. 

In view of the dictionary between polynomials and finite Blaschke products \cite{NgTsang}, it is natural to study Smale's mean value conjecture for finite Blaschke products. Actually, Sheil-Small \cite{sheilsmall} has already pointed out that Smale's proof for $M\le 4$ can be adapted without much modification to show that for finite Blaschke products the same bound $4$ also applies. If we let $K$ and $K_n$ be the counterpart of $M$ and $M_n$ respectively for finite Blaschke products (the precise definition of them will be given in section 2), Sheil-Small also showed that $M \le K$ and $M_n \le K_n$ (see p.365 of \cite{sheilsmall}). Therefore, it would be interesting to consider Smale's mean value conjecture for finite Blaschke products. In this paper, we will prove the following result.

\begin{thm}
$K_n\le2\frac{2n-1+(2n-3)4^{1/(1-n)}}{2n-1}$ for $n\ge2$.
\end{thm}

Unfortunately, $K_n$ is of the form $4 - O({1 \over n})$ as $n \to \infty$ and therefore we cannot get any improvement of $M_n$ for polynomials.

It is known that $M_n=\frac{n-1}{n}$ for $n=2,3$ and $4$ and extremal polynomials exist for each case. The following result shows that this is not true for $K_n$.
  
\begin{thm}
$K_n \ge1$ when $n\ge3$. $K_2=1$ and no extremal finite Blaschke products exist in this case.
\end{thm}

For the dual mean value conjecture, if we let $L$ and $L_n$ be the counterpart of $N$ and $N_n$ respectively for finite Blaschke products (the precise definition will be given in section 2), then we have the following 
  
\begin{thm}
$L_n > 1/4^n$ for $n\ge2$.
\end{thm}

As a corollary, we have an improvement of Dubinin and Sugwa's bound, $N_n \ge \frac{1}{n4^n}$.

\noindent
{\bf Corollary.} $N_n \ge L_n > 1/4^n$ {\it for} $n\ge2$.

Finally, we have  
\begin{thm}
$L_n\le1/n$ when $n\ge3$. $L_2=1/2$ and no extremal finite Blaschke products exist in this case.
\end{thm}

In Section 2, we will review some properties of finite Blaschke products and formulate both Smale's mean value conjecture and the dual mean value conjecture for finite Blaschke products. We will then prove Theorem 1 and 2 in Section 3 and Theorem 3 and 4 in Section 4.

\section{A Brief Review on Finite Blaschke Products}

We recall some basic facts about finite Blaschke products that will be used later. These facts can be found in  \cite{BeardonMinda}, \cite{NgTsang}, \cite{sheilsmall} and the reference therein.
 
\begin{defi}
A \emph{finite Blaschke product} of degree $n$ is a rational function of the form
\begin{equation*}
B(z)=e^{i\alpha}\prod_{k=1}^n \frac{z-z_k}{1-\overline{z_k}z}
\end{equation*}
where $\alpha$ is a real number and $z_1,\ldots,z_n$ are complex numbers on the standard unit disk $\mathbb{D}=\{z:|z|<1\}$.
\end{defi}

It follows immediately from the definition that $B$ has $n$ zeros and $n$ poles, counting multiplicity, and the zeros and poles of $B$ are conjugate relative to the unit circle $\partial\mathbb{D}$. Also we have $|B(z)|<1$ for $|z|<1$, $|B(z)|=1$ for $|z|=1$ and $|B(z)|>1$ for $|z|>1$. If we write $P(z)=\prod_{k=1}^n (z-z_k)$ and define $P^*(z)=z^n\overline{P(1/\overline{z})}$, then $B(z)=e^{i\alpha}\frac{P(z)}{P^*(z)}$. 

Just as polynomials of degree $n$ are precisely the $n$-to-$1$ covering maps of the complex plane onto itself, finite Blaschke products of degree $n$ are precisely the $n$-to-$1$ covering maps of the unit disk $\mathbb{D}$ onto itself. Hence, for each $w\ \in \mathbb{D}$, $B(z)=w$ has exactly $n$ roots in $\mathbb{D}$, counting multiplicity. Let $M(z)=e^{i\theta}\frac{z-\beta}{1-\overline{\beta}z}$ be a M\"obius transformation with $\theta\in\mathbb{R}$ and $\beta\in \mathbb{D}$. Then $B\circ M$ and $M\circ B$ are finite Blaschke products of degree $n$.

By differentiating the relation $B(z)\overline{B(1/\overline{z})}=1$, we conclude that the critical points of $B$ inside and outside the unit disk are conjugating relative to the unit circle. Let 
$\widetilde B (z) = e^{-i\alpha}B(z)$. Then $\widetilde B$ and $B$ share the same set of zeros and critical 
points. Differentiation yields 
$${\widetilde B}'(z)= \frac{P'(z)P^*(z)-P(z)(P^*)'(z)}{\left(P^*(z)\right)^2}=a\frac{Q(z)Q^*(z)}{(P^*(z))^2},$$ 
where $Q(z): =\prod\limits^{n-1}_{i=1} (z- \zeta_i)$ and $\zeta_1, \cdots , \zeta_{n-1}$ 
are the critical points of $\widetilde B$ within the unit disk, and $a$ is a complex constant depending on $B$. This follows from the observation that the numerator is of degree $2n-2$, and the critical points are in conjugating pairs relative to the unit circle, as already remarked. The complex constant $a$ may not be $1$, as one can verify by computing a simple case. Hence we have $B'(z) = e^{i\alpha} aC(z) R^2(z)$
in $\mathbb D$, where $C(z) = \frac{Q(z)}{Q^*(z)}$ 
is a finite Blaschke product of degree $n-1$ with zeros at the critical
points of $B$, and $R(z)= \frac{{Q^*} (z)}{{P^*}(z)}=\frac{\prod_{k=1}^{n-1}(1-\overline{\zeta_k}z)}{\prod_{k=1}^{n-1}(1-\overline{z_k}z)}$ is analytic and non-zero in $\mathbb D$. For an exposition in more details, we refer to Section 11.3 in \cite{sheilsmall}, although there the author omits the complex constant $a$, either for convenience or by error in equation (11.21) of \cite{sheilsmall}.
 
To give a proper formulation of Smale's mean value conjecture for finite Blaschke products in the hyperbolic setting, we introduce the notions of pseudo-hyperbolic distance and hyperbolic derivative.

Let $z$ and $w$ be two points in $\mathbb{D}$. The complex \emph{pseudo-hyperbolic distance} $[z,w]$ is defined by $[z,w]=\frac{z-w}{1-\overline{w}z}$. Let $f$ be a function from $\mathbb{D}$ to $\mathbb{D}$. We define the \emph{hyperbolic derivative} of $f$ at $w$ to be the limit $\lim_{z\rightarrow w} \frac{[f(z),f(w)]}{[z,w]}$, if it exists, and we denote it by $D_H f(w)$.
If $f$ is analytic at a point $w$, then 
$D_H f(w)=\lim_{z\rightarrow w}\frac{\frac{f(z)-f(w)}{1-\overline{f(w)}f(z)}}{\frac{z-w}{1-\overline{w}z}}=\lim_{z\rightarrow w}\frac{f(z)-f(w)}{z-w}\frac{1-\overline{w}z}{1-\overline{f(w)}f(z)}=f'(w)\frac{1-|w|^2}{1-|f(w)|^2}$.

Now we formulate Smale's mean value conjecture for finite Blaschke product in the hyperbolic setting. 
As mentioned before, in the case of complex polynomials, it suffices to consider polynomials with a single zero at $0$ and $a=0$, as we can always apply linear transformations to convert the problem into this case. It happens that we have a similar simplification in the case of finite Blaschke products, provided we formulate the problem in the following way:

\begin{prob}
Let $B$ be a finite Blaschke product of degree $n\ge2$. Suppose $\zeta_1,\ldots, \zeta_{n-1}$ are the critical points of $B$, and suppose $B'(w)\ne0$. For any such $B$, what is the smallest constant $K_n$ such that the following holds for at least one of its critical point $\zeta_i$:
\begin{equation*}
\left|\frac{[B(\zeta_i),B(w)]}{[\zeta_i,w]}\right|\le K_n|D_H B(w)|\text{?}
\end{equation*}
\end{prob}

\begin{rmk}
First of all, it's enough to consider $w=0$. To see this, let $M(z)=\frac{z+w}{1-\overline{w}z}$, and replace $B$ by $C=B\circ M$. Second, we can assume $B(0)=0$ as we can let $M(z)=\frac{z-B(0)}{1-\overline{B(0)}z}$, and replace $B$ by $C=M\circ B$. Finally, we may assume the adjusting coefficient $e^{i\alpha}$ in the expression of $B$ to be $1$. This is easily achieved by replacing $B$ with $e^{-i\alpha}B$.
\end{rmk}

We denote by $\mathcal{B}(n)$ the set of all finite Blaschke products of degree $n$ of the form $B(z)=\prod_{k=1}^n \frac{z-z_k}{1-\overline{z_k}z}$ such that $B(0)=0$ and $B'(0)\ne0$. Then for each $B\in\mathcal{B}(n)$, we define 

$$S(B):=\min\left\{\left|\frac{B(\zeta)}{\zeta B'(0)}\right|:B'(\zeta)=0,\zeta\in \mathbb{D}\right\}.$$

Then by what we have discussed in the previous remark, we have $K_n=\sup\{S(B):B\in\mathcal{B}(n)\}$.

For the dual mean value conjecture for finite Blaschke products, it is clear that we should consider   
$$T(B):=\max\left\{\left|\frac{B(\zeta)}{\zeta B'(0)}\right|:B'(\zeta)=0,\zeta\in \mathbb{D}\right\}$$
and then have $L_n=\inf\{T(B):B\in\mathcal{B}(n)\}$.

\section{Bounds for $K_n$}

The proof of the upper bound for $K_n$ in Theorem 1 follows closely the arguments used by Fujikawa and Sugawa in \cite{fujikawaSugawa} for the bound $M_n\le4\frac{1+(n-2)4^{1/(n-1)}}{n+1}$. The only difference is the algebra involved for finite Blaschke products is more complicated, so we include the detailed arguments here. 

We will need a result in \cite{BeardonMindaNg} concerning the hyperbolic metric for a specific domain. Any simply connected domain $\Omega$ of $\mathbb{C}$ has a hyperbolic metric, which we denote by $\rho_\Omega(z)|dz|$. Let $f$ be a conformal map from $\Omega$ onto another simply connected domain $\Sigma$, then $\rho_\Sigma\left(f(z)\right)|f'(z)|=\rho_\Omega(z)$. In particular we have $\rho_\mathbb{D}(z)=\frac2{(1-|z|)^2}$ for $\mathbb{D}$. Now let $R$ be a positive real number and $l_k$ be the ray of the form $\{re^{i\theta_k}:r\ge R\}$, $k=1,\ldots,n$, where the rays are distinct. For the domain $\Omega:=\mathbb{C}-(l_1\cup\cdots\cup l_n)$, it is proven in \cite{BeardonMindaNg} that $\rho_\Omega(0)\le\frac2R 4^{-1/n}$.


\noindent\emph{Proof of Theorem 1.}\quad Let $B =z\prod_{k=1}^{n-1}\frac{z-z_k}{1-\overline{z_k}z} \in \mathcal{B}(n)$. First note that $B$ is analytic in $\mathbb{C}-\left\{1/\overline{z_1},\ldots,1/\overline{z_{n-1}}\right\}$ and covers every point in $\mathbb{C}$. We may assume the zeros of $B$ are distinct and we know $B$ has exactly $n-1$ critical points in $\mathbb{D}$, counting multiplicity. Let $\zeta_1\ldots\zeta_{n-1}$ be its critical points in $\mathbb{D}$, arranged such that $\min\{|\zeta_k|:k=1,\ldots,n-1\}=|\zeta_1|$. Then the critical points outside $\mathbb{D}$ are precisely $1/\overline{\zeta_1},\ldots,1/\overline{\zeta_{n-1}}$. And we have $B(1/\overline{\zeta_k})=1/\overline{B(\zeta_k)}$. Let $R=\min\{|B(\zeta_k)|:k=1,\ldots,n-1\}>0$. Let $l_k$ be the ray of the form $\{re^{i\theta_k}:r\ge R\}$ that passes through $B(\zeta_k)$ (and hence passes through $B(1/\overline{\zeta_k})$ too). Let $\Omega=\mathbb{C}-(l_1\cup\cdots\cup l_{n-1})$. Consider the inverse branch $\varphi(w)$ of $w=B(z)$ with $\varphi(0)=0$. Then $\varphi(w)$ has expansion $\varphi(w)=a_1w+a_2w+\cdots$ and can be analytically continued to $\Omega$ as it does not contain any critical values of $B$. Hence $\varphi$ is a univalent function in $\Omega$. Let $\mathbb{D}_R$ be the open disk centered at the origin with radius $R$, then $\varphi(w)$ is univalent in $\mathbb{D}_R$. We have $\frac{\varphi(Rw)}{a_1R}$ is a normalized univalent function in $\mathbb{D}$ and $\frac{\varphi(Rw)}{a_1R}=w+\frac{a_2}{a_1}Rw^2+\cdots$. We then have $\left|\frac{a_2}{a_1}R\right|\le2$.

As none of the $\zeta_k$'s lies in $\varphi\left(\mathbb{D}_R\right)$, we have $\varphi_k(w)=\frac{\varphi(w)}{1-\varphi(w)/\zeta_k}$ is univalent throughout $\mathbb{D}_R$. Now $\varphi_k(w)=a_1w+(a_2+a_1^2/\zeta_k)w^2+\cdots$. So we have $\left|\left(\frac{a_2}{a_1}+\frac{a_1}{\zeta_k}\right)R\right|\le2$. Since $\varphi\left(\mathbb{D}_R\right)\subset \mathbb{D}$, we also have $\varphi_x(w)=\frac{\varphi(w)}{1-x\varphi(w)}$ is univalent in $\mathbb{D}_R$ for any $x\in \mathbb{D}$. Hence $\left|\left(\frac{a_2}{a_1}+a_1x\right)R\right|\le2$. By the triangle inequality, $\left|x-1/\zeta_k\right|\le4/(|a_1|R)=4|B'(0)|/R$.

Let $f:\mathbb{D}\rightarrow\Omega$ be a conformal homeomorphism satisfying $f(0)=0$, which has the form $f(u)=c_1u+c_2u^2+\cdots$. The hyperbolic density $\rho_\Omega$ of $\Omega$ at $0$ satisfies $\rho_\Omega(0)|c_1|=\rho_\mathbb{D}(0)$. Hence $|c_1|=\rho_\mathbb{D}(0)/\rho_\Omega(0)\ge R4^{1/(n-1)}$.

Now consider $g(u)=(\varphi\circ f)(u)=a_1c_1u+(a_1c_2+a_2c_1^2)u^2+\cdots$. We have $g(u)$ omits the critical points of $B$. Let $g_x(u)=\frac{g(u)}{1-g(u)/x}=a_1c_1u+(a_1c_2+a_2c_1^2+a_1^2c_1^2/x)u^2+\cdots$, where $x=\zeta_k$ or $1/\overline{\zeta_k}$. Hence $\left|\frac{c_2}{c_1}+c_1\left(\frac{a_2}{a_1}+\frac{a_1}{x}\right)\right|\le2$, therefore $\left|\frac{c_2}{c_1^2}+\frac{a_2}{a_1}+\frac{a_1}{x}\right|\le\frac2{|c_1|}\le\frac2R4^{-1/(n-1)}$. In particular, this holds for $x=\zeta_1$. By the triangle inequality, $\left|\frac1x-\frac1{\zeta_1}\right|\le\frac1{a_1}\left|\left(\frac{c_2}{c_1^2}+\frac{a_2}{a_1}+\frac{a_1}{x}\right)-\left(\frac{c_2}{c_1^2}+\frac{a_2}{a_1}+\frac{a_1}{\zeta_1}\right)\right|\le\frac{|B'(0)|}R4^\frac{n-2}{n-1}$.

By the relation $\varphi\left(B(z)\right)=z$, we have $\varphi''(0)\left(B'(0)\right)^2+\varphi'(0)B''(0)=0$. Hence $2a_2/a_1^2+a_1B''(0)=0$. From Section 2, we know that for $B(z)=z\prod_{k=1}^{n-1}\frac{z-z_k}{1-\overline{z_k}z}$, we have $B'(z)=aC(z)R^2(z)$, where $C(z)=\prod_{k=1}^{n-1}\frac{z-\zeta_k}{1-\overline{\zeta_k}z}$ is the Blaschke product formed by the critical points of $B$, $R(z)=\frac{\prod_{k=1}^{n-1}(1-\overline{\zeta_k}z)}{\prod_{k=1}^{n-1}(1-\overline{z_k}z)}$, and $a$ is a complex constant. It follows that $\frac{1}{a_1}=B'(0)=aC(0)=(-1)^{n-1}a\prod_{k=1}^{n-1}\zeta_k$ and $B''(0)=aC'(0)+2aC(0)R'(0)$.
Now, $R'(z)=\sum_{k=1}^{n-1}\frac{\overline{z_k}-\overline{\zeta_k}}{(1-\overline{z_k}z)^2}\prod_{\substack{1\le s\le n-1\\s\ne k}}\frac{1-\overline{\zeta_s}z}{1-\overline{z_s}z}$ and 
$$C'(z)=\sum_{k=1}^{n-1}\frac{1-\zeta_k\overline{\zeta_k}}{(1-\overline{\zeta_k}z)^2}\prod_{\substack{1\le s\le n-1\\s\ne k}}\frac{z-\zeta_s}{1-\overline{\zeta_s}z}.$$ Hence $R'(0)=\sum_{k=1}^{n-1}(\overline{z_k}-\overline{\zeta_k})$ and
$$C'(0)=(-1)^{n-1}\left(\prod_{k=1}^{n-1}\zeta_k\right)\sum_{k=1}^{n-1}(\overline{\zeta_k}-1/\zeta_k)=C(0)\sum_{k=1}^{n-1}(\overline{\zeta_k}-1/\zeta_k).$$

Summarizing all these, we have $-2a_2/a_1^3=B''(0)=1/a_1\sum_{k=1}^{n-1}(2\overline{z_k}-1/\zeta_k-\overline{\zeta_k})$. Thus

\begin{equation*}
\begin{split}
\frac{4|B'(0)|}R&\ge\left|\frac{2a_2}{a_1^2}+\frac2{\zeta_1}\right|=\left|\frac2{\zeta_1}+\sum_{k=1}^{n-1}\left(-2\overline{z_k}+\frac1{\zeta_k}+\overline{\zeta_k}\right)\right|\\
&=\left|\frac3{\zeta_1}+\sum_{k=1}^{n-1}\left(-2\overline{z_k}+\overline{\zeta_k}\right)+\sum_{k=2}^{n-1}\frac1{\zeta_k}\right|\\
&=\left|\frac{|\zeta_1|}{\zeta_1}\right|\left|\frac3{|\zeta_1|}+\sum_{k=1}^{n-1}\frac{\zeta_1}{|\zeta_1|}\left(-2\overline{z_k}+\overline{\zeta_k}\right)+\sum_{k=2}^{n-1}\frac{\zeta_1}{|\zeta_1|}\frac1{\zeta_k}\right|\\
&\ge \frac3{|\zeta_1|}-2\sum_{k=1}^{n-1}\re{\frac{\zeta_1\overline{z_k}}{|\zeta_1|}}+\sum_{k=1}^{n-1}\re{\frac{\zeta_1\overline{\zeta_k}}{|\zeta_1|}}+\sum_{k=2}^{n-1}\re{\frac{\zeta_1}{|\zeta_1|\zeta_k}}.
\end{split}
\end{equation*}

On the other hand, $4^{\frac{n-2}{n-1}}\frac{|B'(0)|}R\ge\left|\frac1x-\frac1{\zeta_1}\right|\ge\frac1{|\zeta_1|}-\re{\frac{\zeta_1}{|\zeta_1|x}}$, and hence $\re{\frac{\zeta_1}{|\zeta_1|x}}\ge\frac1{|\zeta_1|}-4^{\frac{n-2}{n-1}}\frac{|B'(0)|}R$, where $x=\zeta_k$ or $1/\overline{\zeta_k}$. Also $\frac{4|B'(0)|}R\ge\left|-\overline{z_k}-\frac1{\zeta_k}\right|\ge\frac1{|\zeta_1|}+\re{\frac{\zeta_1\overline{z_k}}{|\zeta_1|}}$ and hence $-\re{\frac{\zeta_1\overline{z_k}}{|\zeta_1|}}\ge\frac1{|\zeta_1|}-4\frac{|B'(0)|}R$. Now we have
\begin{equation*}
\begin{split}
\frac{4|B'(0)|}R&\ge\frac3{|\zeta_1|}-2\sum_{k=1}^{n-1}\re{\frac{\zeta_1\overline{z_k}}{|\zeta_1|}}+\sum_{k=1}^{n-1}\re{\frac{\zeta_1\overline{\zeta_k}}{|\zeta_1|}}+\sum_{k=2}^{n-1}\re{\frac{\zeta_1}{|\zeta_1|\zeta_k}}\\
&\ge\frac3{|\zeta_1|}+2(n-1)\left(\frac1{|\zeta_1|}-\frac{4|B'(0)|}R\right)+(2n-3)\left(\frac1{|\zeta_1|}-4^{\frac{n-2}{n-1}}\frac{|B'(0)|}R\right)\text{,}
\end{split}
\end{equation*}
which gives $\frac1{|\zeta_1|}\le2\frac{2n-1+(2n-3)4^{1/(1-n)}}{2n-1}\frac{|B'(0)|}R$.

\noindent Now for $\zeta_k$ such that $|B(\zeta_k)|=R$, we have
\begin{equation*}
S(B)\le\left|\frac{B(\zeta_k)}{\zeta_kB'(0)}\right|=\frac R{|\zeta_k||B'(0)|}\le\frac R{|\zeta_1||B'(0)|}\le2\frac{2n-1+(2n-3)4^{1/(1-n)}}{2n-1}.
\end{equation*}
\noindent Hence $K_n\le2\frac{2n-1+(2n-3)4^{1/(1-n)}}{2n-1}$ as desired.\hspace*{\fill}$\square$

\bigskip

\noindent\emph{Proof of Theorem 2.}\quad For brevity we write $d=n-1$. Let $\omega_d$ be the primitive $d$-th root of unity. We consider $B(z)=z\prod_{k=0}^{d-1}\frac{z-\omega_d^k \alpha}{1-\overline{\omega}_d^k\alpha z}=z\frac{z^d-\alpha^d}{1-\alpha^dz^d}$, where $\alpha\in(0,1)$. Then $B'(z)=-\frac{\alpha^d z^{2d}+\left((d-1)\alpha^{2d}-(d+1)\right)z^d+\alpha^d}{(1-\alpha^dz^d)^2}$. The critical points of $B$ in the unit disk are $d$-th roots of $\zeta=\frac{(d+1)-(d-1)\alpha^{2d}-\sqrt{\left((d+1)^2-(d-1)^2\alpha^{2d}\right)\left(1-\alpha^{2d}\right)}}{2\alpha^d}=\frac{(d+1)-(d-1)\beta^2-\sqrt{\left((d+1)^2-(d-1)^2\beta^2\right)\left(1-\beta^2\right)}}{2\beta}$, where $\beta=\alpha^d\in(0,1)$. 

Hence, $S(B)=\left|\frac1\beta\frac{\zeta-\beta}{1-\beta\zeta}\right|=\frac1{\beta^2}\frac{\sqrt{(d+1)^2-(d-1)^2\beta^2}-(d+1)\sqrt{1-\beta^2}}{\sqrt{(d+1)^2-(d-1)^2\beta^2}-(d-1)\sqrt{1-\beta^2}}<1$, and tends to $1$ as $\beta$ tends to $1$. Hence $K_n \ge1$. Moreover, it's easy to see that we actually have equality holds when $n=2$ and no extremal finite Blaschke products exist in this case.\hspace*{\fill}$\square$

\section{Bounds of $L_n$}

Dubinin and Sugawa \cite{dubininSugawa} showed that for polynomials, $N_n$ is bounded below by $1/(n4^n)$. Their method can be adapted to prove a similar result for finite Blaschke products (Theorem 3), with the help of the following lemma which improves the inequality (7) in Duren and Schiffer's paper \cite{durenSchiffer}.

\begin{lm}
Let $f$ be a function analytic and univalent in the annulus $A(r)=\{z:r<|z|<1\}$ satisfying the following:
\begin{compactenum}
\item $|f(z)|<1$ in $A(r)$ while $|f(z)|=1$ on $\partial\mathbb{D}$;
\item $f(z)\ne0$ in $A(r)$;
\item $f(1)=1$.
\end{compactenum}
Then $\lim\sup_{|z|\rightarrow r}|f(z)|\le \max\{2r,\frac{4r}{1+4r^2}\}<4r$.
\end{lm}

\noindent\emph{Proof of Lemma 1.}\quad The domain $\mathbb{D}-f(A(r))$ is conformal to the annulus $A(r)$ via the univalent map $f$. Hence the modulus of $\mathbb{D}-f(A(r))$ is given by $\log(1/r)$. Suppose $s=\lim\sup_{|z|\to r}|f(z)|$. Since $f(z)\ne0$ in $A(r)$, $\mathbb{D}-f(A(r))$ is an annular domain separating $0$ and a point with modulus $s$ from the unit circle. Therefore, its modulus is bounded by $\mu(s)$, the modulus of the extreme Gr\"otzsch ring $\mathbb{D}-[0,s]$ (see p.54 of \cite{LV}). From the inequality (4.5) in page 643 of \cite{QV} or the inequality (2.13) in page 63 of \cite{LV}, we have $\mu(s) < \log(\frac{2(1+\sqrt{1-s^2})}{s})$. Hence we have $\log(1/r)<\log(\frac{2(1+\sqrt{1-s^2})}{s})$ or $\frac{s}{r} < 2(1+\sqrt{1-s^2})$. If $s<2r$, then we are done. If $s \ge 2r$, then $(\frac{s}{r}-2)^2 \le 4(1-s^2)$ and hence  $s \le \frac{4r}{1+4r^2}$. Therefore, $\lim\sup_{|z|\rightarrow r}|f(z)|\le \max\{2r,\frac{4r}{1+4r^2}\}<4r$. \hspace*{\fill}$\square$

\noindent\emph{Proof of Theorem 3.}\quad Consider a finite Blaschke product $B=z\prod_{k=1}^{n-1}\frac{z-z_k}{1-\overline{z_k}z} \in \mathcal{B}(n)$ with critical points $\zeta_1,\ldots,\zeta_{n-1}$. Further, assume that $R=|B(\zeta_1)|^{1/n}=\max\{|B(\zeta_k)|^{1/n}:k=1,\ldots,n-1\}$. Consider the annulus $A(R^n)=\{w:R^n<|w|<1\}$. Let $U=\{z:|B(z)|\le R^n\}$, then $U$ contains all the zeros and critical points of $B$. The function $B:\mathbb{D}-U\rightarrow A(R^n)$ is a covering map of degree $n$. We wish to find a bijective holomorphic function $\varphi$ from $A(R)$ to $\mathbb{D}-U$ satisfying the assumptions of Lemma 1. Suppose such a function exists, we have $\lim\sup_{|w|\to R^+}|\varphi(w)| < 4R$ by Lemma 1, and thus $U\subset \mathbb{D}_{4R}$. In particular, $|\zeta_k|< 4R$ and $|z_k| < 4R$. Now $T(B)\ge\left|\frac{B(\zeta_1)}{\zeta_1B'(0)}\right|=\frac{R^n}{|\zeta_1|\prod_{k=1}^{n-1}|z_k|}>\frac1{4^n}$.

It remains to show how to define the function $\varphi$. For this we use the topological theory of covering spaces. The reader may refer to any reference book on algebraic topology for detailed expositions, for example \cite{Hatcher}. Define $\widetilde{A}(R)$ to be the annulus $A(R,1/R)=\{w:R<|w|<1/R\}$ for any $R\in(0,1)$. Now let $R$ assume the same meaning as in the last paragraph. Define $\widetilde{U}:=\{z:|B(z)|\le R^n\text{ or }|B(z)|\ge1/R^n\}$. Then $\widetilde{U}$ contains all the critical points, zeros and poles of $B$ on the complex plane. The restriction of $B$ to $\mathbb{C}-\widetilde{U}$ is a covering map of degree $n$ onto the annulus $\widetilde{A}(R^n)$, sending $1$ to $B(1)$ on the unit circle. The domain $\mathbb{C}-\widetilde{U}$ is connected, as in particular all points in the fiber $B^{-1}(1)$ lie in the same connected component. The induced map $B_*$ is an injective homomorphism from $\pi_1(\mathbb{C}-\widetilde{U},1)$ to $\pi_1(\widetilde{A}(R^n),B(1))\cong\mathbb{Z}$, and the image of this map is a subgroup of index $n$. This subgroup must be $n\mathbb{Z}$.

On the other hand, the map $p_n:\widetilde{A}(R)\to\widetilde{A}(R^n), z\mapsto B(1)z^n$ is also a covering map of degree $n$, sending $1$ to $B(1)$. The induced map on the fundamental group is multiplication by $n$. Hence in particular $(p_n)_*(\pi_1(\widetilde{A}(R),1))=B_*(\pi_1(\mathbb{C}-\widetilde{U},1))$. By the theory of covering spaces (for example, Proposition 1.37 in \cite{Hatcher}), there exists a covering space isomorphism $\varphi:\widetilde{A}(R)\to\mathbb{C}-\widetilde{U}$ such that $\varphi(1)=1$. Since $B\circ\varphi=p_n$ with $B,p_n$ being local homeomorphisms and holomorphic, we conclude that $\varphi$ is holomorphic. It is immediate that $|\varphi(z)|=1$ if $|z|=1$, and $<1$ if $|z|<1$ in its domain of definition.

Restrict the map $\varphi$ to the annulus $A(R)=\{w:R<|w|<1\}$, which we still denote by $\varphi$. This map satisfies all the assumptions of Lemma 1.\hspace*{\fill}$\square$

\noindent
{\it Remark.} The reader may be familiar with the notion of ``smooth covering spaces" introduced in the classical book of Ahlfors and Sario \cite{ahlforsSario}. However, one should be aware that the notion of covering spaces adopted here aligns with the standard notion in algebraic topology, which is different from that in the aforementioned book. Here we require that each point in the base space should possess a neighborhood that is evenly covered. This ensures the homotopy lifting property, which only holds for \emph{regular} smooth covering spaces in the sense of Ahlfors and Sario.

\noindent\emph{Proof of Theorem 4.}\quad Let $C(z)=\frac{z^n-a^n}{1-a^nz^n}$ be a finite Blaschke product of degree $n\ge2$, where $a\in(0,1)$. Let $M(z)=\frac{z+a}{1-az}$. Then $B=C\circ M$ is a finite Blaschke product of degree $n$ with the property $B(0)=0$. Now $C'(z)=\frac{n(1-a^{2n})z^{n-1}}{(1-a^nz^n)^2}$, thus $C(z)$ has one critical point $0$ in $\mathbb{D}$ with multiplicity $n-1$. We have $B'(z)=C'(M(z))M'(z)=C'(M(z))\frac{1-a^2}{(1+az)^2}$, so $B'(z)=0 \iff C'(M(z))=0$. Hence $B(z)$ has one critical point $-a$ in the unit disk. Now $B(-a)=C(M(-a))=C(0)=-a^n$, $B'(0)=C'(M(0))(1-a^2)=C'(a)(1-a^2)=\frac{na^{n-1}}{1-a^{2n}}(1-a^2)$. Hence $\left|\frac{B(-a)}{(-a)B'(0)}\right|=\frac1n\frac{1-a^{2n}}{1-a^2}>\frac1n$, and $\to\frac1n$ as $a\to0$. Therefore $L_n\le1/n$. When $n=2$ it's easily seen that the equality actually holds and no extremal finite Blaschke products exist in this case.\hspace*{\fill}$\square$

\bigskip
\noindent\emph{Proof of the Corollary.} By Theorem 4, it suffices to prove the following lemma. 

\begin{lm} $L_n \le N_n$ for all $n \in \mathbb{N}$.
\end{lm}

\noindent\emph{Proof.} Given $P(z) = z \prod\limits^{n-1}_{i=1} (z-a_i)$ with critical points $c_i$, define $B_m(z) =  z \prod\limits^{n-1}_{i=1} \left (\dfrac{z-\frac{a_i}m} {1-\frac{\overline{a_i}}m z}\right )$. Then when $m$ is sufficiently large, $B_m$ is a finite Blaschke product.
Let $f_m(z) = m^n B_m \left (\dfrac{1}{\ m\ } z\right ) = z \prod\limits^{n-1}_{i=1} \left (\dfrac {z-a_i} {1-\frac{\overline{a_i}}{m^2} z}\right )$. Then $f_m(z) \to P (z)$ locally uniformly on $\mathbb C$ as $m\to \infty$.\\

If $c_{m,i}$ are the critical points of $B_m$, then $d_{m,i} = m c_{m,i}$ are the critical points of $f_m$ and 
$$  \dfrac{f_m(d_{m,i})}{d_{m,i} f'_m (0)} = 
\dfrac{m^n B_m (\frac 1{\ m \ } d_{m,i})}{d_{m,i} m^{n-1} B'_m (0)}= \dfrac{B_m (c_{m,i})}{c_{m,i} B'_m (0)}.$$

So if for some $C>0$ such that there exists some $1\le j\le n-1$ such that 
$$\left | \dfrac{B_m (c_{m,j})}{c_{m,j} B'_m (0)} \right | \ge C,$$\\
then as $m\to \infty$, we also have there exists some $1\le k \le n-1$ such that
$$\left | \dfrac{ P (c_k)}{c_k P' (0)} \right | \ge C.$$

By the definition of $N_n$ and $L_n$, we must have $N_n \ge L_n$.\hspace*{\fill}$\square$

\noindent
{\it Remark.} The connection between polynomials and finite Blaschke products through some rescalings was pointed to the first author by Toshiyuki Sugawa and Yum Tong Siu, independently.   

\noindent
{\it Remark.} There is a much more direct way to get $N_n \ge 1/4^n$ (but not $N_n > 1/4^n$). The idea is that inequality (2.3) of \cite{dubininSugawa} remains true if one replaces the critical points in it by zeros of the polynomials.   

\section{Conclusion}

Smale's mean value conjecture for finite Blaschke products seems to be as difficult as the original Smale's conjecture. The algebra involved can be quite complicated for finite Blaschke products. Moreover, neither the derivative of a finite Blaschke product nor the sum of two finite Blaschke products is a finite Blaschke products. So some arguments for polynomials cannot be carried over to finite Blaschke products. However, one may hope that    
there are more results from geometric function theory that one can apply as finite Blaschke products are self maps of the unit disk.
For example, we have the following  

\begin{prop}
Let $B$ be a finite Blaschke product of degree $n$ with $B(0) =0$, $B'(0) \ne 0$. If $r = \min\limits_{B'(\zeta_i) =0}  |B (\zeta_i) | \ge \dfrac 1{\,2\,}$, then there
exists some critical points $\zeta_k$ and $\zeta_l$ such that
$$\left | \dfrac{ B(\zeta_k)}{\zeta_k B'(0)} \right | \le 4^{\frac{n-2}{n-1}}\dfrac  2 {\ 3 \ }$$
and 
$$\left | \dfrac{ B(\zeta_l)}{\zeta_l B'(0)} \right | \ge \dfrac 1 {\ 2^n \ }\, ,$$
where $n$ is the degree of the finite Blaschke product $B$.
\end{prop}

\noindent
{\it Proof}. From the proof of Theorem 1, we can extract the following lemma easily.

\begin{lm}
Let $B$ be a finite Blaschke product of degree $n$ with $B(0) =0$, $B'(0) \ne 0$ 
and $r = \min\limits_{B'(\zeta_i) =0}  |B (\zeta_i) | = |B(\zeta_k ) | >0$ 
for some $1\le k \le n-1 $. Let $\varphi = B^{-1} :\mathbb D_r \to \mathbb D$ be such that $\varphi(0)=0$.
If $x,y \in \mathbb C$ are such that $(1-x \varphi (w)) (1- y \varphi (w) )\ne 0$
for all $w\in \mathbb D_r$, then
$$|x-y| \le 4^{\frac{n-2}{n-1}}\dfrac{|B'(0)|} r.$$
\end{lm}

Now take  $ x = \dfrac 1 {\zeta_k}$, $y = \dfrac{e^{i\theta} B (\zeta_k)}{\zeta_k}$ where 
$\theta$ is chosen such that  
$e^{i\theta} B (\zeta_k) = - | B(\zeta_k)|$.
Note that $| y | <1$ as $\dfrac{B(z)}z$ is a finite
Blaschke product.

By Lemma 3, we have 
$$\left | \dfrac{ 1+ |B(\zeta_k)|}{\zeta_k} \right | \le 4^{\frac{n-2}{n-1}}\dfrac {|B'(0)|} r.$$
Hence
$$\left | \dfrac{B(\zeta_k)}{\zeta_k B'(0)} \right | \le \dfrac{4^{\frac{n-2}{n-1}}}{1+|B(\zeta_k)|} \, $$
and the first inequality follows.

For the second inequality, since $r = \min\limits_{B'(\zeta_i) =0}  |B (\zeta_i) | \ge \dfrac 1{\,2\,}$, $\max\limits_{B'(\zeta_i) =0}  |B (\zeta_i) | \ge \dfrac 1{\,2^n\,}$ and hence $R:=\max\limits_{B'(\zeta_i) =0}  |B (\zeta_i) |^{1/n}= |B (\zeta_l) |^{1/n} \ge \dfrac 1{\,2\,}$. It follows that $2R \ge \frac{4R}{1+4R^2}$. As in the proof of Theorem 3, applying Lemma 1, we have $|\zeta_k|\le2R$ and $|z_k|\le2R$ for all $1 \le k \le n-1$. Now $\left|\frac{B(\zeta_l)}{\zeta_lB'(0)}\right|=\frac{R^n}{|\zeta_l|\prod_{k=1}^{n-1}|z_k|}\ge\frac1{2^n}$ and we are done. \hspace*{\fill}$\square$

\medskip
\noindent
{\it Acknowledgment.} The authors want to thank the referee for her/his many very helpful suggestions, in particular, a very detailed suggestion of a more topological proof of Theorem 3. We have followed the suggestion of making use of the topological theory of covering spaces as well as applying a better estimate of $\mu(s)$, the modulus of the extreme Gr\"otzsch ring $\mathbb{D}-[0,s]$ in Section 4. We also thank Toshiyuki Sugawa and Yum Tong Siu for the suggestion of the connection between polynomials and finite Blaschke products through some rescalings.

\end{document}